\def\version{0.16}

\def\journal{XXX}
\def\titlep{A tensor product of representations
of UHF algebras arising from Kronecker products
}
\documentclass[11pt]{article}
\usepackage{graphicx,ifthen}
\usepackage{amssymb}
\usepackage{amsmath}

\font\germ=eufm10 at12pt

\def\goth#1{\hbox{\germ#1}}

\newcommand{\inlim}{{\displaystyle \lim_{\to}\,}}

\newcommand{\qed}{\hbox{\rule[-2pt]{3pt}{6pt}}}
\newcommand{\qedh}{\hfill\qed \\}

\newcommand{\vv}{\vspace{.3in}}

\newtheorem{Thm}{Theorem}[section]

\newtheorem{rem}[Thm]{Remark}

\newtheorem{prop}[Thm]{Proposition}
\newtheorem{prob}[Thm]{Problem}
\newtheorem{cor}[Thm]{Corollary}
\newtheorem{fact}[Thm]{Fact}

\newcommand{\kn}{\Large\bf
$K\hspace{-.4cm} N$
\Large\bf\vv }

\def\cal#1{\mathcal #1}
\def\con{{\cal O}_{n}}

\def\pr{{\it Proof.}\quad}

\def\co#1{{\cal O}_{#1}}
\def\brl{branching law}
\def\bfsnl{{\rm BFS}_{N}(\Lambda)}

\addtocounter{footnote}{1}
\def\sftt#1{
\setcounter{equation}{0}
\addtocounter{footnote}{1}
\section{#1}
}

\def\ssft#1{\subsection{#1}}

\def\sssft#1{\subsubsection{#1}}

\begin{document}
%
%
\def\autherp{Katsunori Kawamura}
\def\emailp{e-mail: kawamura@kurims.kyoto-u.ac.jp.}
\def\addressp{{\small {\it College of Science and Engineering, 
Ritsumeikan University,}}\\
{\small {\it 1-1-1 Noji Higashi, Kusatsu, Shiga 525-8577, Japan}}
}

\def\infw{\Lambda^{\frac{\infty}{2}}V}
\def\zhalfs{{\bf Z}+\frac{1}{2}}
\def\ems{\emptyset}
\def\pmvac{|{\rm vac}\!\!>\!\! _{\pm}}
\def\vac{|{\rm vac}\rangle _{+}}
\def\dvac{|{\rm vac}\rangle _{-}}
\def\ovac{|0\rangle}
\def\tovac{|\tilde{0}\rangle}
\def\expt#1{\langle #1\rangle}
\def\zph{{\bf Z}_{+/2}}
\def\zmh{{\bf Z}_{-/2}}
\def\brl{branching law}
\def\bfsnl{{\rm BFS}_{N}(\Lambda)}
\def\scm#1{S({\bf C}^{N})^{\otimes #1}}
\def\mqb{\{(M_{i},q_{i},B_{i})\}_{i=1}^{N}}
\def\zhalf{\mbox{${\bf Z}+\frac{1}{2}$}}
\def\zmha{\mbox{${\bf Z}_{\leq 0}-\frac{1}{2}$}}
\newcommand{\mline}{\noindent
\thicklines
\setlength{\unitlength}{.1mm}
\begin{picture}(1000,5)
\put(0,0){\line(1,0){1250}}
\end{picture}
\par
 }
\def\ptimes{\otimes_{\varphi}}
\def\delp{\Delta_{\varphi}}
\def\delps{\Delta_{\varphi^{*}}}
\def\gamp{\Gamma_{\varphi}}
\def\gamps{\Gamma_{\varphi^{*}}}
\def\sem{{\sf M}}
\def\hdelp{\hat{\Delta}_{\varphi}}
\def\tilco#1{\tilde{\co{#1}}}
\def\ndm#1{{\bf M}_{#1}(\{0,1\})}
\def\cdm#1{{\cal M}_{#1}(\{0,1\})}
\def\tndm#1{\tilde{{\bf M}}_{#1}(\{0,1\})}
\def\sck{{\sf CK}_{*}}
\def\hdel{\hat{\Delta}}
\def\ba{\mbox{\boldmath$a$}}
\def\bb{\mbox{\boldmath$b$}}
\def\bc{\mbox{\boldmath$c$}}
\def\bd{\mbox{\boldmath$d$}}
\def\be{\mbox{\boldmath$e$}}
\def\bk{\mbox{\boldmath$k$}}
\def\bm{\mbox{\boldmath$m$}}
\def\bp{\mbox{\boldmath$p$}}
\def\bq{\mbox{\boldmath$q$}}
\def\bu{\mbox{\boldmath$u$}}
\def\bv{\mbox{\boldmath$v$}}
\def\bw{\mbox{\boldmath$w$}}
\def\bx{\mbox{\boldmath$x$}}
\def\by{\mbox{\boldmath$y$}}
\def\bz{\mbox{\boldmath$z$}}
\def\bomega{\mbox{\boldmath$\omega$}}
\def\N{{\bf N}}
\def\lxm{L_{2}(X,\mu)}
\def\titlepage{

\noindent
{\bf 
\noindent
\thicklines
\setlength{\unitlength}{.1mm}
\begin{picture}(1000,0)(0,-300)
\put(0,0){\kn \knn\, for \journal\, Ver.\version}
\put(0,-50){\today}
\end{picture}
}
\vspace{-2.3cm}
\quad\\
{\small file: \textsf{tit01.txt,\, J1.tex}
 \footnote{
  ${\displaystyle
  \mbox{directory: \textsf{\fileplace}, 
  file: \textsf{\incfile},\, from \startdate}}$
          }
}
\quad\\
\framebox{
 \begin{tabular}{ll}
 \textsf{Title:} &
 \begin{minipage}[t]{4in}
 \titlep
 \end{minipage}
 \\
 \textsf{Author:} &\autherp
 \end{tabular}
}
{\footnotesize	\tableofcontents }
}
\def\ngt{{\bf N}_{\geq 2}}
\def\ngti{\ngt^{\infty}}
%
%
%
\setcounter{section}{0}
\setcounter{footnote}{0}
\setcounter{page}{1}
\pagestyle{plain}

%
%
\title{\titlep}
\author{\autherp\thanks{\emailp}
\\
\addressp}
\date{}
\maketitle

%
%
\begin{abstract}
We introduce a non-symmetric tensor product
of representations of UHF algebras
by using Kronecker products of matrices.
We prove tensor product formulae of  GNS representations by product states
and show examples.
\end{abstract}


%
%
\sftt{Introduction}
\label{section:first}
We have studied representations of operator algebras.
We found non-symmetric tensor products of representations of some 
C$^{*}$-algebras 
which are associated with
non-cocommutative comultiplication 
of certain C$^{*}$-bialgebras \cite{TS01,TS05,TS07}.
In this paper, 
we introduce a tensor product of representations of 
uniformly hyperfinite (=UHF) algebras by using
Kronecker products of matrices.
According to the tensor product,
we show tensor product formulae of 
Gel'fand-Na\u{\i}mark-Segal (=GNS) representations by product states
with respect to given factorizations of UHF algebras.
In this section, we show our motivation,
definitions and main theorems.

%
%
\ssft{Motivation}
\label{subsection:firstone}
In this subsection, we roughly explain our motivation 
and the background of this study.
Explicit mathematical definitions will 
be shown after $\S$ \ref{subsection:firsttwo}.

For any group $G$,
there always exists the standard inner tensor product 
(or Kronecker product) of representations of $G$ 
associated with the diagonal map from $G$ to
$G\times G$ \cite{FH,Pukanszky}. The tensor product
is important to describe the duality of $G$ \cite{Tatsuuma}.
For an algebra $A$,
one does not know how to define 
an associative inner tensor product of representations 
of $A$ unless $A$ has a coassociative comultiplication.
In \cite{TS01}, we introduced a non-symmetric tensor product 
among all representations of Cuntz algebras and determined tensor product 
formulae of all permutative representations completely,
in spite of the unknown of any comultiplication of Cuntz algebras.
In \cite{TS02}, we generalized this construction of 
tensor product to a system of C$^{*}$-algebras
and $*$-homomorphisms indexed by a monoid.
For example, we constructed a non-symmetric tensor product
of all representations of Cuntz-Krieger algebras 
by using Kronecker products of matrices \cite{TS11,TS05,TS07}.
About the basic idea of tensor products of
representations of C$^{*}$-algebras, see $\S$ 1.1 in \cite{TS01}.

On the other hand,
UHF algebras and their representations 
are well studied \cite{AK02R,AW,A,AN,BJ,Glimm, Powers,N}.
For example,
GNS representations of product states of UHF algebras 
were completely classified by \cite{AN}.
This class contains representations of UHF algebras
of all Murray-von Neumann's types I, II, III (\cite{Blackadar2006}, $\S$ III.5).

Our interest is to construct a tensor product of representations of UHF algebras 
and compute tensor product formulae of this class of representations.
Since one knows neither a cocommutative nor non-commutative comultiplication 
of UHF algebra, the tensor product is new if one can find it.

%
%
\ssft{GNS representations of UHF algebras by product states}
\label{subsection:firsttwo}
We recall GNS representations by product states of UHF algebras
which are given as infinite tensor products of matrix algebras
\cite{AW,A,N}.
Let $\ngt\equiv \{2,3,4,\ldots\}$ and 
let $\ngti$ denote the set of all sequences in $\ngt$.
For $1\leq n<\infty$,
let $M_{n}$ denote the (finite-dimensional)
C$^{*}$-algebra of all $n\times n$-complex matrices.
For $\ba= (a_{n})_{n\geq 1}\in \ngti$,
the sequence $\{M_{a_{n}}\}_{n\geq 1}$ of C$^{*}$-algebras
defines the tensor product 
%
%
\begin{equation}
\label{eqn:finite}
{\cal A}_{n}(\ba)\equiv \bigotimes_{j=1}^{n}M_{a_{j}}.
\end{equation}
With respect to the embedding 
%
%
\begin{equation}
\label{eqn:psi}
\psi_{\ba}^{(n)}:{\cal A}_{n}(\ba)\ni A\mapsto
 A\otimes I\in {\cal A}_{n}(\ba)\otimes M_{a_{n+1}}
={\cal A}_{n+1}(\ba),
\end{equation}
we regard ${\cal A}_{n}(\ba)$ as a C$^{*}$-subalgebra of 
${\cal A}_{n+1}(\ba)$ and
let 
${\cal A}(\ba)$ denote the inductive limit of the system
$\{({\cal A}_{n}(\ba),\psi_{\ba}^{(n)}):n\geq 1\}$:
%
%
\begin{equation}
\label{eqn:ab}
{\cal A}(\ba)\equiv  \inlim ({\cal A}_{n}(\ba),\psi_{\ba}^{(n)}).
\end{equation}
By definition, ${\cal A}(\ba)$ is a UHF algebra of 
Glimm's type $\{a_{1}\cdots a_{k}\}_{k\geq 1}$ 
which was classified by \cite{Glimm}.
On the contrary, any UHF algebra is isomorphic to ${\cal A}(\ba)$ 
for some $\ba\in \ngti$.
Hence we call ${\cal A}(\ba)$ a UHF algebra in this paper.

Let $M_{n,+,1}$ denote the set of all positive elements in $M_{n}$ 
which traces are $1$.
Then any state of $M_{n}$ is written as 
$\omega_{T}(x)\equiv {\rm tr}(Tx)$ ($x\in M_{n}$)
for some $T\in M_{n,+,1}$
where ${\rm tr}$ denotes the trace of $M_{n}$.
Let 
${\cal T}(\ba)\equiv \prod_{n\geq 1}M_{a_{n},+,1}$ and 
let $\{E^{(n)}_{ij}:i,j=1,\ldots,n\}$ denote the set 
of standard matrix units of $M_{n}$.
For ${\bf T}=(T^{(n)})_{n\geq 1}\in {\cal T}(\ba)$,
define the state
$\omega_{{\bf T}}$ of ${\cal A}(\ba)$ by
%
%
\begin{equation}
\label{eqn:tstate}
\omega_{{\bf T}}(E_{j_{1}k_{1}}^{(a_{1})}
\otimes\cdots\otimes E_{j_{n}k_{n}}^{(a_{n})}
)\equiv T_{k_{1},j_1}^{(1)}\cdots T_{k_{n},j_n}^{(n)}
\end{equation}
for each $j_{1},\ldots,j_{n},k_{1},\ldots,k_{n}$ and $n\geq 1$
where $T^{(n)}_{jk}$'s denote matrix elements of the matrix $T^{(n)}$.
Then $\omega_{{\bf T}}$ coincides with the product state 
$\bigotimes_{n\geq 1}\omega_{T^{(n)}}$.
%
%
\begin{Thm}
\label{Thm:awone}(\cite{K-R}, Remark 11.4.16)
For each ${\bf T}\in {\cal T}(\ba)$,
the state $\omega_{{\bf T}}$ in (\ref{eqn:tstate})
is a factor state, that is,
if $({\cal H}_{{\bf T}},\pi_{{\bf T}},x_{{\bf T}})$ is 
the GNS triplet of ${\cal A}(\ba)$
by $\omega_{{\bf T}}$,
then $\pi_{{\bf T}}({\cal A}(\ba))^{''}$ is a factor.
\end{Thm}
The factor ${\cal M}_{{\bf T}}\equiv \pi_{{\bf T}}({\cal A}(\ba))^{''}$ 
is called an {\it Araki-Woods factor (or infinite tensor product 
of finite dimensional type {\rm I} (=ITPFI) factor)} \cite{AW, A}.
Properties of ${\cal M}_{{\bf T}}$ 
and $({\cal H}_{{\bf T}},\pi_{{\bf T}},x_{{\bf T}})$ 
are closely studied in \cite{AW,A,Shlyakhtenko} and \cite{AN}, respectively.

%
%
\ssft{A set of isomorphisms}
\label{subsection:firstthree}
In this subsection, we introduce a set of isomorphisms among algebras
${\cal A}(\ba)$ in (\ref{eqn:ab}) and their tensor products.
By using the set,
we will define a tensor product of representations and that of states
in $\S$ \ref{subsection:firstfour}.

For $\ba=(a_{n}),\bb=(b_{n})\in \ngti$,
let $\ba\cdot \bb\equiv (a_{1}b_{1},a_{2}b_{2},\ldots)\in\ngti$ and
define the $*$-isomorphism $\varphi_{\ba,\bb}^{(n)}$ 
of ${\cal A}_{n}(\ba\cdot \bb)$
onto ${\cal A}_{n}(\ba)\otimes {\cal A}_{n}(\bb)$ by
%
%
\begin{equation}
\label{eqn:embeddingtwo}
\varphi_{\ba,\bb}^{(n)}(E_{j_{1}k_{1}}^{(a_{1}b_{1})}\otimes \cdots \otimes
E_{j_{n}k_{n}}^{(a_{n}b_{n})})
\equiv 
(E_{j_{1}^{'},k_{1}^{'}}^{(a_{1})}\otimes\cdots\otimes 
E_{j_{n}^{'},k_{n}^{'}}^{(a_{n})})
\otimes 
(E_{j_{1}^{''},k_{1}^{''}}^{(b_{1})}\otimes\cdots\otimes 
E_{j_{n}^{''},k_{n}^{''}}^{(b_{n})})
\end{equation}
for each $j_{i},k_{i}\in \{1,\ldots,a_{i}b_{i}\}$, $i=1,\ldots,n$
where
$j_{1}^{'},\ldots,j_{n}^{'}$,
$k_{1}^{'},\ldots,k_{n}^{'}$,
$j_{1}^{''},\ldots,j_{n}^{''}$,
$k_{1}^{''},\ldots,k_{n}^{''}$
are defined as
$j_{i}=b_{i}(j_{i}^{'}-1)+j_{i}^{''}$, 
$k_{i}=b_{i}(k_{i}^{'}-1)+k_{i}^{''}$
and
$j_{i}^{'},k^{'}_{i}\in \{1,\ldots,a_{i}\}$,
$j_{i}^{''},k_{i}^{''}\in \{1,\ldots,b_{i}\}$
for each $i=1,\ldots,n$.
For $\psi_{\ba}^{(n)}$ in (\ref{eqn:psi}), we see that
%
%
\begin{equation}
\label{eqn:psitwo}
(\psi_{\ba}^{(n)}\otimes \psi_{\bb}^{(n)})\circ 
\varphi^{(n)}_{\ba,\bb}
=
\varphi^{(n+1)}_{\ba,\bb}\circ \psi_{\ba\cdot \bb}^{(n)}
\end{equation}
for each $\ba,\bb$ and $n$.
From $\{\varphi_{\ba}^{(n)}:n\geq 1\}$,
we can define the $*$-isomorphism $\varphi_{\ba,\bb}$ 
of ${\cal A}(\ba\cdot \bb)$
onto ${\cal A}(\ba)\otimes {\cal A}(\bb)$ such that 
%
%
\begin{equation}
\label{eqn:isomorphism}
(\varphi_{\ba,\bb})|_{{\cal A}_{n}(\ba\cdot \bb)}=\varphi_{\ba,\bb}^{(n)}
\end{equation}
for each $n$
where we identify ${\cal A}(\ba)\otimes {\cal A}(\bb)$ 
with the inductive limit of the system
$\{({\cal A}_{n}(\ba)\otimes {\cal A}_{n}(\bb),\psi^{(n)}_{\ba}\otimes \psi^{(n)}_{\bb}):n\geq 1\}$.
Then the following holds:
%
%
\begin{equation}
\label{eqn:commutative}
(\varphi_{\ba,\bb}\otimes id_{\bc})\circ \varphi_{\ba\cdot \bb,\bc}
=
(id_{\ba}\otimes \varphi_{\bb,\bc})\circ \varphi_{\ba,\bb\cdot \bc}
\quad(\ba,\bb,\bc\in \ngti)
\end{equation}
where $id_{\bx}$ denotes the identity map on ${\cal A}(\bx)$
for $\bx=\ba,\bc$.
Equivalently,
the following diagram is commutative:

\noindent
%
%
%
\thicklines
\setlength{\unitlength}{.1mm}
\begin{picture}(1000,300)(-30,120)
\put(-30,250){${\cal A}(\ba\cdot \bb\cdot \bc)$}
\put(150,280){\vector(3,1){200}}
\put(150,350){$\varphi_{\ba,\bb\cdot \bc}$}
\put(150,240){\vector(3,-1){200}}
\put(150,160){$\varphi_{\ba\cdot \bb,\bc}$}
\put(390,350){${\cal A}(\ba)\otimes {\cal A}(\bb\cdot \bc)$}
\put(390,150){${\cal A}(\ba\cdot \bb)\otimes {\cal A}(\bc)$}
\put(650,350){\vector(3,-1){200}}
\put(710,350){$id_{\ba}\otimes \varphi_{\bb,\bc}$}
\put(650,170){\vector(3,1){200}}
\put(710,160){$\varphi_{\ba,\bb}\otimes id_{\bc}$}
\put(880,250){${\cal A}(\ba)\otimes {\cal A}(\bb)\otimes {\cal A}(\bc)$.}
\end{picture}
%
%
\begin{rem}
\label{rem:one}
{\rm
\begin{enumerate}
\item
We consider the meaning of (\ref{eqn:embeddingtwo}).
For two matrices $A\in M_{n}$ and $B\in M_{m}$,
define the matrix $A\boxtimes B\in M_{nm}$
by 
%
%
\begin{equation}
\label{eqn:box}
(A\boxtimes B)_{m(i-1)+i^{'},m(j-1)+j^{'}}\equiv 
A_{ii^{'}}B_{jj^{'}}
\end{equation}
for $i,i^{'}\in \{1,\ldots,n\}$
and $j,j^{'}\in \{1,\ldots,m\}$.
The new matrix $A\boxtimes B$ 
is called the {\it Kronecker product} of $A$ and $B$ \cite{Dief}.
When $n=1$, we see that
%
%
\begin{equation}
\label{eqn:oneone}
\varphi_{\ba,\bb}^{(1)}(A\boxtimes B)
=A \otimes B\quad(A\in M_{a_{1}},\,B\in M_{b_{1}}).
\end{equation}
Hence $\varphi_{\ba,\bb}^{(1)}$ is the inverse operation of the Kronecker product,
which should be called the {\it Kronecker coproduct}.
\item
It is clear that ${\cal A}(\ba\cdot \bb)$ and 
${\cal A}(\ba)\otimes {\cal A}(\bb)$ are isomorphic
even if one does not know $\varphi_{\ba,\bb}$.
However, the choice of isomorphisms $\{\varphi_{\ba,\bb}\}$ 
in (\ref{eqn:isomorphism}),
which satisfies both (\ref{eqn:psitwo}) 
and (\ref{eqn:commutative}) is not trivial.
\end{enumerate}
}
\end{rem}

%
%
\ssft{Main theorems}
\label{subsection:firstfour}
In this subsection, we show our main theorems.
%
%
\sssft{Basic properties}
\label{subsubsection:firstfourone}
In this subsubsection, 
we introduce a tensor product
of representations 
and that of states of UHF algebras and 
show its basic properties.
For a C$^{*}$-algebra ${\goth A}$,
let ${\rm Rep}{\goth A}$
and ${\cal S}({\goth A})$ denote
the class of all representations
and the set of all states of ${\goth A}$, respectively.
By using the set $\{\varphi_{\ba,\bb}:\ba,\bb\in \ngti\}$
in (\ref{eqn:isomorphism}),
define the operation $\ptimes$ from 
${\rm Rep}{\cal A}(\ba)\times {\rm Rep}{\cal A}(\bb)$
to ${\rm Rep}{\cal A}(\ba\cdot \bb)$ by
%
%
\begin{equation}
\label{eqn:ptimes}
\pi_{1}\ptimes \pi_{2}\equiv (\pi_{1}\otimes \pi_{2})\circ \varphi_{\ba,\bb}
\end{equation}
for $(\pi_{1},\pi_2)\in {\rm Rep}{\cal A}(\ba)\times {\rm Rep}{\cal A}(\bb)$.
We see that 
if $\pi_{i}$ and $\pi_{i}^{'}$
are unitarily equivalent for $i=1,2$,
then 
$\pi_{1}\ptimes \pi_{2}$
and 
$\pi_{1}^{'}\ptimes \pi_{2}^{'}$
are also unitarily equivalent.
Furthermore, 
define the operation $\ptimes$ from 
${\cal S}({\cal A}(\ba))\times {\cal S}({\cal A}(\bb))$
to ${\cal S}({\cal A}(\ba\cdot \bb))$ by
%
%
\begin{equation}
\label{eqn:rho}
\rho_{1}\ptimes \rho_{2}\equiv (\rho_{1}\otimes \rho_{2})\circ \varphi_{\ba,\bb}
\end{equation}
for $(\rho_{1},\rho_2)\in 
{\cal S}({\cal A}(\ba))\times {\cal S}({\cal A}(\bb))$.
From (\ref{eqn:commutative}),
we see that
%
%
\begin{equation}
\label{eqn:ptimestwo}
(\pi_{1}\ptimes \pi_{2})\ptimes \pi_{3}=
\pi_{1}\ptimes (\pi_{2}\ptimes \pi_{3}),\quad
(\rho_{1}\ptimes \rho_{2})\ptimes \rho_{3}=
\rho_{1}\ptimes (\rho_{2}\ptimes \rho_{3})
\end{equation}
for each 
$(\pi_{1},\pi_2,\pi_{3})\in {\rm Rep}{\cal A}(\ba)\times {\rm Rep}{\cal A}(\bb)
\times {\rm Rep}{\cal A}(\bc)$
and 
for $(\rho_{1},\rho_2,\rho_{3})\in
{\cal S}({\cal A}(\ba))\times {\cal S}({\cal A}(\bb))
\times {\cal S}({\cal A}(\bc))$
and $\ba,\bb,\bc\in \ngti$.

The following fact 
is a paraphrase of well-known results of tensor products of factors
(which will be proved in $\S$ \ref{section:second}).
%
%
\begin{fact}
\label{fact:types}
Let $\pi_{1}$ and $\pi_{2}$
be representations of 
${\cal A}(\ba)$
and 
${\cal A}(\bb)$, respectively.
Then the following holds:
\begin{enumerate}
\item
If both $\pi_{1}$ and $\pi_{2}$
are factor representations,
then so is $\pi_{1}\ptimes \pi_{2}$.
\item
The type of $\pi_{1}\ptimes \pi_{2}$
coincides with 
that of $\pi_{1}\otimes \pi_{2}$
where the type of a representation $\pi$ of a C$^{*}$-algebra ${\goth A}$
means the type of the von Neumann algebra $\pi({\goth A})^{''}$
(\cite{Blackadar2006}, Theorem III.2.5.27).
\item
If both $\pi_{1}$ and $\pi_{2}$
are irreducible,
then 
so is $\pi_{1}\ptimes \pi_{2}$.
\end{enumerate}
\end{fact}

By definition,
the essential part of the tensor product $\ptimes$
is given by the set $\{\varphi_{\ba,\bb}\}$ of isomorphisms
in (\ref{eqn:embeddingtwo}).
The idea of the definition of $\{\varphi_{\ba,\bb}\}$
is an analogy of the set of embeddings of Cuntz algebras
in $\S$ 1.2 of \cite{TS01}.
This type of tensor product is known yet
in neither operator algebras nor the purely algebraic theory 
of quantum groups \cite{Kassel}.
%
%
\begin{rem}
\label{rem:second}
{\rm
\begin{enumerate}
\item
Our terminology ``tensor product of representations" is different
from usual sense \cite{FH}.
Remark that,
for $\pi,\pi^{'}\in {\rm Rep}{\cal A}(\ba)$,
$\pi\ptimes \pi^{'}\not\in {\rm Rep}{\cal A}(\ba)$ but
$\pi\ptimes \pi^{'}\in {\rm Rep}{\cal A}(\ba\cdot \ba)$
because $\ba\cdot \ba\ne \ba$ for any $\ba\in \ngti$.
\item
From Fact \ref{fact:types}(iii),
there is no nontrivial branch of 
the irreducible decomposition of 
the tensor product of any two irreducibles.
In general, such a tensor product of the other algebra is decomposed
into more than one irreducible component.
For example, see Theorem 1.6 of \cite{TS01}.
\end{enumerate}
}
\end{rem}

%
%
\sssft{Tensor product formulae of GNS representations}
\label{subsubsection:firstfourtwo}
Next, we show tensor product formulae of GNS representations
of product states in Theorem \ref{Thm:awone} as follows.
%
%
\begin{Thm}
\label{Thm:main}
Let  $\ba,\bb\in \ngti$ and 
let $\omega_{{\bf T}}$ be as in (\ref{eqn:tstate})
with the GNS representation $\pi_{{\bf T}}$.
\begin{enumerate}
\item
For each ${\bf T}\in {\cal T}(\ba)$ and ${\bf R}\in {\cal T}(\bb)$, 
%
%
\begin{equation}
\label{eqn:omega}
\omega_{{\bf T}}\ptimes \omega_{{\bf R}}=\omega_{{\bf T}\boxtimes {\bf R}}
\end{equation}
where 
${\bf T}\boxtimes {\bf R}\in {\cal T}(\ba\cdot \bb)$ is
defined as
%
%
\begin{equation}
\label{eqn:boxt}
{\bf T}\boxtimes {\bf R}\equiv (T^{(1)}\boxtimes R^{(1)},
T^{(2)}\boxtimes R^{(2)},T^{(3)}\boxtimes R^{(3)},\ldots)
\end{equation}
for ${\bf T}=(T^{(n)})$ and ${\bf R}=(R^{(n)})$.
\item
For each ${\bf T}=(T^{(n)})\in {\cal T}(\ba)$ 
and ${\bf R}=(R^{(n)})\in {\cal T}(\bb)$, 
$\pi_{{\bf T}}\ptimes \pi_{{\bf R}}$ 
is unitarily equivalent
to $\pi_{{\bf T}\boxtimes {\bf R}}$.
\end{enumerate}
\end{Thm}

\noindent
From Theorem \ref{Thm:main},
the tensor product $\ptimes$ is compatible with 
product states and their GNS representations.
More precisely, for the following 
two semigroups $({\cal T},\boxtimes)$ and $({\cal S},\ptimes)$, the map 
%
%
\begin{equation}
\label{eqn:calt}
{\cal T}\equiv \bigcup_{\ba\in \ngti}{\cal T}(\ba)
\ni {\bf T}\mapsto \omega_{{\bf T}}\in {\cal S}\equiv 
\bigcup_{\ba\in \ngti}{\cal S}({\cal A}({\ba}))
\end{equation}
is a semigroup homomorphism.
Let ${\sf R}_{\ba}$ denote the set of all unitary equivalence classes
in ${\rm Rep}{\cal A}(\ba)$.
Then 
%
%
\begin{equation}
\label{eqn:caltt}
{\cal T}\ni {\bf T}\mapsto [\pi_{{\bf T}}]\in {\sf R}\equiv 
\bigcup_{\ba\in \ngti}{\sf R}_{{\ba}}
\end{equation}
is also a  semigroup homomorphism between
$({\cal T},\boxtimes)$ and $({\sf R},\ptimes)$
where $[\pi]$ denotes the unitary equivalence class of 
a representation $\pi$.

From Theorem \ref{Thm:main}, the following holds.
%
%
\begin{cor}
\label{cor:nonsymmetric}
\begin{enumerate}
\item
There exist 
$\ba,\bb\in\ngti$, and states
$\omega\in {\cal S}({\cal A}(\ba))$ and
$\omega^{'}\in {\cal S}({\cal A}(\bb))$ such that
$\omega\ptimes \omega^{'}\ne \omega^{'}\ptimes \omega$.
\item
There exist 
$\ba,\bb\in\ngti$, and
classes
$[\pi]\in {\sf R}_{\ba}$ and
$[\pi^{'}]\in {\sf R}_{\bb}$
such that
$[\pi]\ptimes [\pi^{'}]\ne 
[\pi^{'}]\ptimes [\pi]$.
\end{enumerate}
\end{cor}

\noindent
From Corollary \ref{cor:nonsymmetric}(i) and (ii),
we say that $\ptimes$ is non-symmetric.
In other words,
two semigroups $({\cal S},\ptimes)$
and $({\sf R},\ptimes)$ are non-commutative.
These non-commutativities come from the
non-commutativity of the Kronecker product of matrices.

%
%
\begin{prob}
\label{prob:af}
{\rm
\begin{enumerate}
\item
Generalize the tensor product in (\ref{eqn:ptimes})
to that of representations of AF algebras
which are not always UHF algebras.
\item
Reconstruct UHF algebras from 
$({\cal S},\ptimes)$ and $({\sf R},\ptimes)$,
and show a Tatsuuma duality type theorem 
for UHF algebras \cite{Tatsuuma}.
\end{enumerate}
}
\end{prob}

In $\S$ \ref{section:second},
we will prove theorems in $\S$ \ref{subsection:firstfour}.
In $\S$ \ref{section:third}, we will show examples of Theorem \ref{Thm:main}.

%
%
\sftt{Proofs of main theorems}
\label{section:second}
In this section, we prove main theorems.
\\
\\
{\it Proof of Fact \ref{fact:types}.}
By the definition of $\ptimes$,
%
%
\begin{equation}
\label{eqn:productone}
(\pi_{1}\ptimes \pi_{2})({\cal A}(\ba\cdot \bb))
=
(\pi_{1}\otimes \pi_{2})({\cal A}(\ba)\otimes {\cal A}(\bb)).
\end{equation}
\noindent
(i)
Since $\pi_{1}\otimes \pi_{2}$ is also a factor representation,
the statement holds from (\ref{eqn:productone}).

\noindent
(ii)
By definition,
the type of $\pi_{1}\ptimes \pi_{2}$
is the type of 
$\{(\pi_{1}\ptimes \pi_{2})({\cal A}(\ba\cdot \bb))\}^{''}$.
From this and (\ref{eqn:productone}), the statement holds.

\noindent
(iii) By assumption, $\pi_{1}\otimes \pi_{2}$ is also irreducible.
From this and (\ref{eqn:productone}), the statement holds.
\qedh

\noindent
{\it Proof of Theorem \ref{Thm:main}.}
(i)
By Definition \ref{eqn:tstate},
the statement holds from direct computation.

\noindent
(ii)
Let $\ba\in \ngti$.
For ${\bf T}\in {\cal T}(\ba)$,
let $({\cal H}_{{\bf T}},\pi_{{\bf T}},\Omega_{{\bf T}})$
denote the GNS triplet by the state $\omega_{{\bf T}}$.
Define the GNS map $\Lambda_{{\bf T}}$  \cite{KV,MNW}
from ${\cal A}(\ba)$ to ${\cal H}_{{\bf T}}$ by
%
%
\begin{equation}
\label{eqn:lambda}
\Lambda_{{\bf T}}(x)\equiv \pi_{{\bf T}}(x)\Omega_{{\bf T}}
\quad(x\in {\cal A}({\bf T})).
\end{equation}
Let $\ba,\bb\in \ngti$.
For ${\bf T}\in {\cal T}(\ba)$ and ${\bf R}\in {\cal T}(\bb)$,
define the unitary $U^{({\bf T},{\bf R})}$
from ${\cal H}_{{\bf T}\boxtimes {\bf R}}$
to ${\cal H}_{{\bf T}}\otimes {\cal H}_{{\bf R}}$ by
%
%
\begin{equation}
\label{eqn:utr}
U^{({\bf T},{\bf R})}
\Lambda_{{\bf T}\boxtimes {\bf R}}(x)
\equiv 
(\Lambda_{{\bf T}}
\otimes
\Lambda_{{\bf R}})(\varphi_{\ba,\bb}(x))
\quad(x\in {\cal A}(\ba\cdot \bb)).
\end{equation}
Since $\varphi_{\ba,\bb}$ is bijective,
$U^{({\bf T},{\bf R})}$ is well-defined as a unitary,
and we see that
%
%
\begin{equation}
\label{eqn:ur}
U^{({\bf T},{\bf R})}\pi_{{\bf T}\boxtimes {\bf R}}(x)
(U^{({\bf T},{\bf R})})^{*}
=
(\pi_{{\bf T}}\otimes \pi_{{\bf R}})(\varphi_{\ba,\bb}(x))
=
(\pi_{{\bf T}}\ptimes \pi_{{\bf R}})(x)
\end{equation}
for $x\in {\cal A}(\ba\cdot \bb)$.
Hence two representations
$\pi_{{\bf T}\boxtimes {\bf R}}$ and 
$\pi_{{\bf T}}\ptimes \pi_{{\bf R}}$
are unitarily equivalent.
\qedh

\noindent
{\it Proof of Corollary \ref{cor:nonsymmetric}.}
Let $\ba=\bb=(2,2,2,\ldots)\in\ngti$.
Define $T,R\in M_{2,+,1}$ by
$T\equiv {\rm diag}(1,0)$ and
$R\equiv {\rm diag}(0,1)$.
Define
${\bf T}=(T^{(n)})\in {\cal T}(\ba)$ and
${\bf R}=(R^{(n)})\in {\cal T}(\bb)$ 
by  $T^{(n)}=T$ and $R^{(n)}=R$ 
for all $n\geq 1$.
Then 
$T^{(n)}\boxtimes R^{(n)}=T\boxtimes R={\rm diag}(0,1,0,0)$
and 
$R^{(n)}\boxtimes T^{(n)}=R\boxtimes T={\rm diag}(0,0,1,0)$.
Let  $\omega\equiv \omega_{{\bf T}}$
and $\omega^{'}\equiv \omega_{{\bf R}}$.
Then
%
%
\begin{equation}
\label{eqn:minus}
(\omega\ptimes \omega^{'})(E_{22}^{(4)}-E_{33}^{(4)})=1
\ne -1=(\omega^{'}\ptimes \omega)(E_{22}^{(4)}-E_{33}^{(4)})
\end{equation}
where $\{E^{(4)}_{ij}\}$ denotes the set of standard matrix units of $M_{4}$.
This implies (i).

From (\ref{eqn:minus}), 
$\|\omega\ptimes \omega^{'}-\omega^{'}\ptimes \omega\|=2$.
From  Corollary 10.3.6 of \cite{K-R},
GNS representations $\pi_{0}$ and $\pi^{'}_{0}$
of ${\cal A}(\ba\cdot \bb)$ by
$\omega\ptimes \omega^{'}$ and 
$\omega^{'}\ptimes \omega$ are disjoint.
Especially, 
$\pi_{0}$ and $\pi_{0}^{'}$ are not unitarily equivalent.
From Theorem \ref{Thm:main}(i) and (ii),
$\pi_{0}$ and $\pi^{'}_{0}$ are unitarily equivalent to
$\pi_{{\bf T}}\ptimes \pi_{{\bf R}}$
and 
$\pi_{{\bf R}}\ptimes \pi_{{\bf T}}$,
respectively.
Therefore
%
%
\begin{equation}
\label{eqn:pit}
[\pi_{{\bf T}}]\ptimes [\pi_{{\bf R}}]=[\pi_{0}]
\ne [\pi^{'}_{0}]=[\pi_{{\bf R}}]\ptimes [\pi_{{\bf T}}].
\end{equation}
Hence (ii) holds.
\qedh

%
%
\sftt{Example}
\label{section:third}
In this section,
we show an example of Theorem \ref{Thm:main}.
Let ${\cal A}(\ba)$ and ${\cal T}(\ba)$ be as in $\S$ \ref{subsection:firsttwo}.
For $n\geq 2$,
let $\ba(n)\equiv (n,n,n,\ldots)\in \ngti$ and
let 
%
%
\begin{equation}
\label{eqn:uhfn}
UHF_{n}\equiv {\cal A}(\ba(n))=(M_{n})^{\otimes \infty}.
\end{equation}
Then $UHF_{n}$ is 
the UHF algebra of Glimm's type $\{n^{l}\}_{l\geq 1}$.
We consider tensor product formulae of a small class of representations of $UHF_{n}$.

Let $\{E^{(n)}_{ij}\}$ denote the set of all standard matrix units 
of $M_{n}$.
For $j\in \{1,\ldots,n\}$,
define $F_{j}^{(n)}\equiv E_{jj}^{(n)}$.
Then $F_{j}^{(n)}\in M_{n,+,1}$ for each $j$.
For $J=(j_{l})_{l\geq 1}\in \{1,\ldots,n\}^{\infty}$,
define 
${\bf T}(J)\equiv (F_{j_{1}}^{(n)},F_{j_{2}}^{(n)},\ldots)\in {\cal T}(\ba(n))$.
%
%
\begin{prop}
\label{prop:subclass}
Let $\pi_{{\bf T}}$ be as in Theorem \ref{Thm:awone}.
Then the following holds:
\begin{enumerate}
\item
For each $J\in \{1,\ldots,n\}^{\infty}$,
$\pi_{{\bf T}(J)}$ is irreducible.
\item
We write $P_{n}[J]$ as the unitary equivalence class of $\pi_{{\bf T}(J)}$.
Then
%
%
\begin{equation}
\label{eqn:productb}
P_{n}[J]\ptimes P_{m}[K]=P_{nm}[J\cdot K]
\end{equation}
for each 
$J=(j_{l})\in \{1,\ldots,n\}^{\infty}$,
$K=(k_{l})\in \{1,\ldots,m\}^{\infty}$
and $n,m\geq 2$ 
where
$J\cdot K\in {\cal T}(\ba(nm))$ is defined by
%
%
\begin{equation}
\label{eqn:dotdot}
J\cdot K\equiv (m(j_{1}-1)+k_{1},m(j_{2}-1)+k_{2},m(j_{3}-1)+k_{3},\ldots).
\end{equation}
\end{enumerate}
\end{prop}
%
%
\pr
(i)
Since the state $\omega_{{\bf T}(J)}$ is the product state
of pure states,
$\omega_{{\bf T}(J)}$ is pure.
Hence the statement holds.

\noindent
(ii)
Recall $\boxtimes$ in (\ref{eqn:boxt}). 
Then we can verify that
${\bf T}(J)\boxtimes {\bf T}(K)={\bf T}(J\cdot K)$
for each
$J\in \{1,\ldots,n\}^{\infty}$ and
$K\in \{1,\ldots,m\}^{\infty}$.
From this and Theorem \ref{Thm:main}(ii),
the statement holds.
\qedh

\noindent
From Proposition \ref{prop:subclass},
${\sf P}\equiv \bigcup_{n\geq 2}\{P_{n}[J]:J\in\{1,\ldots,n\}^{\infty}\}$
is a semigroup 
of unitary equivalence classes of irreducible representations
with respect to the product $\ptimes$. 
From the proof of Corollary \ref{cor:nonsymmetric},
$({\sf P},\ptimes)$ is also non-commutative.

The class $P_{n}[J]$ in Proposition \ref{prop:subclass}(ii)
coincides with
the restriction of a permutative representation of the Cuntz algebra $\con$
on the UHF subalgebra of $U(1)$-gauge invariant elements in $\con$, 
which is called an atom \cite{AK02R,BJ}.
This class contains only type I representations of $UHF_{n}$, but
relations with representations of Cuntz algebras 
and quantum field theory are well studied \cite{AK02R,BJ}.


%
%

\end{document}